\documentstyle[12pt]{article}
\newtheorem{theorem}{Theorem}[section]

\newcommand{\U}{U_q(C_n^{(1)})}
\renewcommand{\a}{\alpha}

\newcommand{\la}{\lambda}
\newcommand{\be}{\beta}

\newcommand{\ep}{\epsilon}
\newcommand{\vep}{\varepsilon}

\title{Bosonic Realizations of $U_q(C^{(1)}_n)$}

\author{Naihuan Jing, Yoshitaka Koyama and Kailash C. Misra}
\begin{document}           % End of preamble and beginning of text.
%\date{July 11, 1996}
\maketitle                 % Produces the title.

\begin{center}
Department of Mathematics\\
North Carolina State University\\
Raleigh, NC 27695-8205, USA
\end{center}

\begin{abstract}
We construct explicitly the quantum symplectic affine algebra
$U_q(\widehat{sp}_{2n})$ using bosonic fields. The Fock
space decomposes into irreducible modules of level $-1/2$, quantizing
the Feingold-Frenkel construction for $q=1$.
\end{abstract}

\section{Introduction}
Quantum affine algebras play an important role in understanding
the trigonometric solutions of the quantum Yang-Baxter equation \cite{kn:Dr}.
Subsequent work shows that they also serve as
symmetry for many statistical models \cite{kn:JM}. In most of the
theoretical and applicational tasks one needs to
construct them
explicitly by known infinite dimensional algebras.

The first such explicit bosonic construction was done for all the
simply laced $(ADE)^{(1)}$ untwisted types by Frenkel and Jing
\cite{kn:FJ}, and subsequently the twisted types $(ADE)^{(r)}$ \cite{kn:J1}
and type $B_n^{(1)}$ \cite{kn:Br} were also completed. The fermionic
constructions were furnished in \cite{kn:H}. The $q$-Wakimoto construction
was also known \cite{kn:M, kn:ABE, kn:S} afterwards. However, the quantum
symplectic affine algebra $U_q(C_n^{(1)})$ was left without an
explicit construction except for $n=1$ \cite{kn:K}.

In the classical $(q=1)$ case Kac-Wakimoto \cite{kn:KW1, kn:KW2}
have introduced an wider class of irreducible highest
weight representations of affine Lie algebras called admissible
representations. Unlike
integrable highest weight representations these representations admit
fractional levels.
Feingold and Frenkel \cite{kn:FF} constructed level $-1/2$
highest weight admissible representations of the symplectic
affine Lie algebra.

In this paper we will construct explicit realizations of
the quantum affine algebra $U_q(C_n^{(1)})$ at level $-1/2$.
Our realization can be thought as a
$q$-analog of Feingold-Frenkel construction,
but it appeals to the equivalent
form of homogeneous quantum Z-algebra construction considered in
\cite{kn:J2}. These are achieved
in terms of the internal bosonic fields and auxiliary bosonic
fields unlike the classical $\beta\gamma$-system, which is inappropriate for
quantization.
Due to this new phenomenon we use
certain screening operators to show that
the Fock space contains four irreducible highest weight
modules of level $-1/2$. These are $q$-analog of
the admissible representations considered by Feingold-Frenkel in
\cite{kn:FF}, where they used a different method to get the irreducibility.

\section{Quantum affine algebras $U_q(C^{(1)}_n)$}

Let
 $A=(A_{ij}), i,j\in I=\{0, 1, \cdots, n\}$ be the  generalized Cartan matrix
of type
$C_n^{(1)}$ so that
\begin{equation}\label{E:2.1}
A=2\sum_{i=0}^nE_{ii}-\sum_{i=0}^{n-1}(E_{i,i+1}+E_{i+1,i})-E_{10}-
E_{n-1, n},
\end{equation}
where $E_{ij}$'s are the unit matrices in ${\bf Z}^{(n+1)\times (n+1)}$.

Let $\a_i=\vep_i-\vep_{i+1}$ ($i=1, \cdots, n-1$) and $\a_n=2\vep_n$
be the simple roots of the simple Lie algebra ${sp}_{2n}$, and $
\lambda_i=\vep_1+\cdots+\vep_i$ ($i=1, \cdots, n$) be the fundamental
weights. Let $P={\bf Z}\vep_1+\cdots+{\bf Z}\vep_n$ and
$Q={\bf Z}\a_1+\cdots+{\bf Z}\a_n$ be the weight and root lattices.
We then let $\Lambda_i, i=0, \cdots n$ be the fundamental
weights for the affine Lie algebra $\widehat{sp}_{2n}$, here $\Lambda_i=
\lambda_i+\Lambda_0$.
The nondegenerate symmetric bilinear form $(\ |\ )$ on ${\bf h}^*$, the
daul Cartan subalgebra of $\widehat{sp}_{2n}$,
satisfies that
\begin{equation}
(\vep_i|\vep_j)=\frac12\delta_{ij},\ \
(\alpha_i|\alpha_j)=d_iA_{ij}, \ \ (\delta|\alpha_i)=(\delta|\delta)=0
\ \ \mbox{for all} \ i,j\in I \label{E:2.2},
\end{equation}
where $(d_0, \cdots, d_n)=(1, 1/2, \cdots, 1/2, 1)$.

Let $q_i=q^{d_i}=q^{\frac 12(\a_i|\a_i)}, i\in I$.
We now recall Drinfeld's realization of the quantum affine
algebra $U_q(C_n^{(1)})$ \cite{kn:Dr}.
The quantum affine algebra $U_q(C_n^{(1)})$ is
the associative algebra with 1 over ${\bf C}(q^{1/2})$
generated by the elements $x_{ik}^{\pm}$, $a_{il}$, $K_i^{\pm 1}$,
$\gamma^{\pm 1/2}$, $q^{\pm d}$ $(i=1,2,\cdots,n, k\in {\bf Z},
l\in {\bf Z} \setminus \{0\})$ with the following defining relations :
\begin{eqnarray}
 [\gamma^{\pm 1/2}, u]&=&0 \ \ \mbox{for all} \ u\in U_q(C_n^{(1)}),
\label{E:R1}\\
\mbox{} [a_{ik}, a_{jl}]&=&\delta_{k+l,0}
\displaystyle\frac {[A_{ij}k]_i}{k}
\displaystyle\frac {\gamma^k-\gamma^{-k}}{q_j-q^{-1}_j},
\label{E:R2}\\
\mbox{} [a_{ik}, K_j^{\pm 1}]&=&[q^{\pm d}, K_j^{\pm 1}]=0,
\label{E:R3}\\
 q^d x_{ik}^{\pm} q^{-d}&=&q^k x_{ik}^{\pm }, \ \
q^d a_{il} q^{-d}=q^l a_{il},
\label{E:R4}\\
 K_i x_{jk}^{\pm} K_i^{-1}&=&q^{\pm (\alpha_i|\alpha_j)} x_{jk} ^{\pm},
\label{E:R5}\\
\mbox{} [a_{ik}, x_{jl}^{\pm}]&=&\pm \displaystyle\frac {[A_{ij}k]_i}{k}
\gamma^{\mp |k|/2} x_{j,k+l}^{\pm},
\label{E:R6}\\
(z-q^{\pm(\a_i|\a_j)}w)X_i^{\pm}(z)X_j^{\pm}(w)
           &+&(w-q^{\pm(\a_i|\a_j)}z)X_j^{\pm}(w)X_i^{\pm}(z) =0
\label{E:R7}\\
\mbox{}
[X_i^+(z),X_j^-(w)]
             =\frac{\delta_{ij}}{(q_i-q_i^{-1})zw}
                &&
                 \left(
                 \psi_i(w\gamma^{1/2})
                 \delta(\frac{w\gamma}{z})
                 -\varphi_i(w\gamma^{-1/2})
                  \delta(\frac{w\gamma^{-1}}z)
                 \right)
\label{E:R8}
\end{eqnarray}
where
$X_i^{\pm}(z)=\sum_{n\in {\bf Z}}x^{\pm}_{i,n}z^{-n-1}$, $\psi_{im}$ and
$\varphi_{im}$ $(m\in {\bf Z}_{\ge 0})$
are defined by
\begin{eqnarray}
&&\psi_i(z)=\sum_{m=0}^{\infty} \psi_{im} z^{-m}
=K_i \textstyle {exp} \left( (q_i-q_i^{-1}) \sum_{k=1}^{\infty} a_{ik}
z^{-k}
\right),
\nonumber\\
&&\phi_i(z)=\sum_{m=0}^{\infty} \varphi_{i,-m} z^{m}
=K_i^{-1} \textstyle {exp} \left(- (q_i-q_i^{-1}) \sum_{k=1}^{\infty}
a_{i,-k} z^{k}
\right),
\nonumber\\
&&\mbox{Sym}_{k_1,\cdots,k_m}\sum_{r=0}^{m=1-A_{ij}}(-1)^r
\left[\begin{array}{c} m\\r\end{array}\right]_iX_i^{\pm}(z_1)\cdots
X^{\pm}_{i}(z_r) \nonumber\\
&&\qquad\qquad\qquad \times X^{\pm}_{j}(w)X^{\pm}_{i}(z_{r+1})\cdots
X^{\pm}_{i}
(z_m)=0
\label{E:R9}
\end{eqnarray}

\section{Fock representation of $U_q(C^{(1)}_n)$}

Let $a_i(m)$ be the operators satisfying the Heisenberg relations for $\U$ at
$\gamma=q^{-1/2}$ and $b_i(m)$ be the independent free bosonic operators:
\begin{equation}
\begin{array}{rcl}
\ [a_i(m), a_j(l)]&=&\delta_{m+l, 0}\frac{[mA_{ij}]_i}{m}
\frac{q^{-m/2}-q^{m/2}}{q_j-q_j^{-1}},
\quad [b_i(m), b_j(l)]=m\delta_{ij}\delta_{m+l,0},\\
\ [a_i(m), b_j(l)]&=&0 .
\end{array}
\end{equation}

   We define the module ${\cal F}_{\a,\be}={\cal F}^1_{\a}\otimes
{\cal F}^2_{\be}$ for the Heisenburg algebra
   by the defining relations
   $$ a_i(m) | \a, \be \rangle = 0 \quad (m>0)
      \; , \quad
      b_i(m) | \a, \be \rangle = 0 \quad (m>0) \;, $$
   $$ a_i(0) | \a, \be \rangle = (\a_i|\a) | \a, \be \rangle
      \; , \quad
      b_i(0) | \a, \be \rangle = (2\vep_i|\be) | \a, \be \rangle \; , $$
   where $| \a, \be \rangle=| \a\rangle\otimes |\be \rangle$
 $(\a\in  P+\frac{{\bf Z}}2\la_n, \be\in P)$ is
   the vacuum vector.
   The grading operator ${d}$ is defined by
 $$ d. | \alpha,\beta \rangle
      =( (\alpha|\alpha)-(\beta|\beta-\lambda_n))
       |\alpha,\beta \rangle . $$

   We set the Fock space
   $$ \widetilde{\cal F}:=
          \bigoplus_{\a\in P+\frac12{\bf Z}\la_n,\be \in P}
          {\cal F}_{\a, \be} . $$
   Let $e^{a_i}$ and $e^{b_i}$ be operators
   on ${\widetilde{\cal F}}$ given by:
   $$ e^{a_i}|\a,\be \rangle= |\a+\a_i,\be \rangle
      \quad , \quad
      e^{b_i}|\a, \be \rangle= |\a,\be+\vep_i \rangle . $$
   Let $: \quad :$ be the usual bosonic normal ordering defined by
   $$ :a_i(m) a_j(l): = a_i(m) a_j(l) \, (m \leq l) ,
   \; a_j(l) a_i(m) \, (m>l) , $$
   $$ :e^{a} a_i(0):=:a_i(0) e^{a}:=e^{a} a_i(0) \, , $$
   $$ :e^{b} b_i(0):=:b_i(0) e^{b}:=e^{b} b_i(0) \, . $$
   and similar normal products for the $b_i(m)$'s.

Let $\partial=\partial_{q^{1/2}}$ be the $q$-difference operator:
$$\partial_{q^{1/2}}(f(z))=\frac{f(q^{1/2}z)-f(q^{-1/2}z)}{(q^{1/2}-q^{-1/2})z}
$$
   We introduce the following operators.
\begin{eqnarray*}
&Y_i^{\pm}(z)&=
         \exp ( \pm \sum^{\infty}_{k=1}
                 \frac{a_i(k)}{[-\frac{1}{2d_i}k]_i} q^{\pm \frac{k}{4}} z^k)
         \exp ( \mp \sum^{\infty}_{k=1}
           \frac{a_i(k)   }{[-\frac{1}{2d_i}k]_i} q^{\pm \frac{k}{4}} z^{-k})
         e^{\pm a_i} z^{\mp 2a_i(0)} ,\\
&Z^{\pm}_i(z)&=
             \exp ( \pm \sum^{\infty}_{k=1} \frac{b_i(-k)}{k} z^k)
             \exp ( \mp \sum^{\infty}_{k=1} \frac{b_i(k)}{k} z^{-k})
             e^{\pm b_i} z^{\pm b_i(0)} .
\end{eqnarray*}
\medskip
\begin{theorem} \label{T:1}The Fock space
    $\widetilde{\cal F}$ is a $U_q$-module of level $-\frac{1}{2}$
    under the action defined by
    $ \gamma \longmapsto q^{-\frac{1}{2}},
       K_i \longmapsto q^{a_i(0)}$,
    $ a_{im} \longmapsto a_i(m),
       q^d \longmapsto q^{{d}}$, and
\begin{eqnarray*}
X_i^+(z) &\longmapsto&
             \partial Z^+_i(z)Z^-_{i+1}(z)Y_i^+(z), \qquad i=1, \cdots n-1\\
X_i^-(z) &\longmapsto&
             Z^-_i(z)\partial Z^+_{i+1}(z)Y^-_i(z), \qquad i=1, \cdots n-1\\
X_n^+(z) &\longmapsto&
             \left(
             \frac{1}{q^{\frac{1}{2}}+ q^{-\frac{1}{2}}}
             :Z^+_n(z) \partial^2_z Z^+_n(z):
     -:\partial_z Y^+_b(q^{\frac{1}{2}}z)
              \partial_z Y^+_b(q^{-\frac{1}{2}}z):
             \right)
             Y^+_a(z) \\
X_n^-(z) &\longmapsto&
             \frac{1}
                  {q^{\frac{1}{2}}+ q^{-\frac{1}{2}}}
             :Z^-_n(q^{\frac{1}{2}}z)Z^-_n(q^{-\frac{1}{2}}z): Y^-_n(z).
\end{eqnarray*}
\end{theorem}

We now prove the theorem by checking that the defined action satisfy
the Drinfeld relations. It is easy to see that the relation
(\ref{E:R1}-\ref{E:R5}) are true by the construction. The relation
(\ref{E:R6}) follows from the expression of $Y^{\pm}_i(z)$ and the
commutativity
between $Y^{\pm}_i(z)$ and $Z_j^{\pm}(z)$. We are left to show the
relations (\ref{E:R7}-\ref{E:R9})

To this end we first list the relations for $Y^{\pm}_i(z)$ and $Z^{\pm}(z)$:
$$
\begin{array}{rcl}
&&Y_i^{\pm}(z)Y_j^{\pm}(w)=:Y_i^{\pm}(z)Y_j^{\pm}(w):\\
   &&\quad\times  \left\{\begin{array}{ll}
     1,                &\mbox{if $(\a_i|\a_j)=0$}\\
     (z-q^{\pm 1/2}w), &\mbox{if $(\a_i|\a_j)=-1/2$}\\
     ((z-q^{\pm 1}w)(z-w))^{-(\a_i|\a_j)},
                             &\mbox{if $ (\a_i|\a_j)=\pm 1$}\\
     \left((z-w)(z-qw)(z-q^{-1}w)(z-q^{\pm 2}w)\right)^{-1},
                         &\mbox{if $ (\a_i|\a_j)=2$}
     \end{array}
     \right. \\
     &&\\
&&Y_i^{\pm}(z)Y_j^{\mp}(w)
=:Y_i^{\pm}(z)Y_j^{\mp}(w):\\
   &&\quad\times  \left\{\begin{array}{ll}
     1,                &\mbox{if $(\a_i|\a_j)=0$}\\
     (z-w)^{-1}, &\mbox{if $(\a_i|\a_j)=-1/2$}\\
     \left((z-q^{-1/2}w)(z-q^{1/2}w)\right)^{(\a_i|\a_j)},
                             &\mbox{if $ (\a_i|\a_j)=\pm 1$}\\
     (z-q^{-1/2}w)(z-q^{1/2}w)(z-q^{-3/2}w)(z-q^{3/2}w),
                         &\mbox{if $ (\a_i|\a_j)=2$}
     \end{array}
     \right.\\
&&\\    && Z_i^{\ep}(z)Z_i^{\ep'}(w)
=:Z_i^{\ep}(z)Z_j^{\ep'}(w):(z-w)^{\ep\ep'\delta_{ij}}.
\end{array}
$$
where the contraction factors such as
$(z-w)^{-1}$ are understood as power series in $w/z$.

For $\epsilon=\pm=\pm 1$ we define
$$\begin{array}{rcl}
X^{+}_{i\epsilon}(z)&=&Z^+_i(q^{\epsilon/2}z)Z^-_{i+1}(z)Y_i^+(z)\\
X^{-}_{i\epsilon}(z)&=&Z^-_i(z)Z^+_{i+1}(q^{\epsilon/2}z)Y_i^-(z)\\
X^{+}_{n\epsilon}(z)&=&:Z^+_n(q_{\epsilon}^{1+\epsilon}z)
Z^+_{n}(q_{\epsilon}^{-1+\epsilon}z)Y_n^+(z):,
\end{array}
$$
where $i=1, \cdots, n-1$ and we allow $\epsilon=0$ in $X^{+}_{n\epsilon}(z)$,
and $q_{\pm}=q^{1/2}$, $q_0=q$. Then we have
$$\begin{array}{rcl}
X_i^{\pm}(z)&=&\frac1{(q^{1/2}-q^{-1/2})z}
\left(X_{i+}^{\pm}(z)-X_{i+}^{\pm}(z)\right)\\
&&\\
X_n^{+}(z)&=&\frac{-1}{(q-q^{-1})(q^{1/2}-q^{-1/2})z^2}
\left(q^{1/2}X_{n+}^{+}(z)+q^{-1/2}X_{n-}^{+}(z)-(q^{1/2}+q^{-1/2})X_{n0}^+(z)
\right)
\end{array}
$$
where $i=1, \cdots, n-1$.

By Wick's theorem we can easily derive operator product expansions for
$X^{\pm}_i(z)$. For instance
we have
\begin{eqnarray}\label{E:R11}
X^{\pm}_{i\ep}(z)X^{\pm}_{i\ep'}(w)
&=&:X^{\pm}_{i\ep}(z)X^{\pm}_{i\ep'}(w):
      (q^{\ep/2}z-q^{\ep'/2}w)(z-q^{\pm 1}w)^{-1}\\
X^{+}_{i\ep}(z)X^{-}_{i\ep'}(w)
&=&:X^{+}_{i\ep}(z)X^{-}_{i\ep'}(w):
      (q^{\ep/2}z-w)^{-1}(z-q^{-\ep'/2}w) \label{E:R12}\\
X^{-}_{i\ep}(z)X^{+}_{i\ep'}(w)
&=&:X^{-}_{i\ep}(z)X^{+}_{i\ep'}(w):
      (z-q^{-\ep/2}w)(q^{\ep'/2}z-w)^{-1} \label{E:R13}\\
X^{+}_{i\ep}(z)X^{+}_{i+1,\ep'}(w)
&=&:X^{+}_{i\ep}(z)X^{+}_{i+1,\ep'}(w):
      (z-q^{\ep'/2}w)^{-1}(z-q^{1/2}w) \label{E:R14}\\
X^{-}_{i\ep}(z)X^{-}_{i+1,\ep'}(w)
&=&:X^{-}_{i\ep}(z)X^{-}_{i+1,\ep'}(w):
      (q^{\ep/2}z-w)^{-1}(z-q^{-1/2}w) \label{E:R15}\\
X^{+}_{i\ep}(z)X^{-}_{i+1,\ep'}(w)&=&:X^{+}_{i\ep}(z)X^{-}_{i+1,\ep'}(w):
\nonumber\\
X^{+}_{i+1,\ep}(z)X^{-}_{i\ep'}(w)
&=&:X^{+}_{i+1,\ep}(z)X^{-}_{i\ep'}(w):
      (q^{\ep/2}z-q^{\ep'/2}w)(z-w)^{-1} \label{E:R16}\\
X^{+}_{n\ep}(z)X^-_{n}(w)&=&:X^{+}_{n\ep}(z)X^-_{n}(w):
\left(\frac{q^{-3\ep/2}(z-q^{3/2}w)}{q^{-\ep/2}z-w}\right)^{|\ep|}
\label{E:R18}
\end{eqnarray}
where one can also
read $X^{+}_{i+1,\ep}(z)X^{+}_{i\ep'}(w)=:X^{+}_{i+1,\ep}(z)X^{+}_{i\ep'}(w):
      (z-q^{\ep/2}w)^{-1}(z-q^{1/2}w)$ etc.

{\it Proof of (\ref{E:R7}).} Most of $+$ cases and $-$ cases are similar.
We show the representative cases in the following.
\begin{eqnarray*}
&&(z-q^{-1/2}w)X^+_{i\ep}(z)X^+_{i+1,\ep'}(w)\\
&=&:X^+_{i\ep}(z)X^+_{i+1,\ep'}(w):
(z-q^{\ep'/2}w)^{-1}(z-q^{-1/2}w)(z-q^{1/2}w)\\
&=&(q^{-1/2}z-w)X^+_{i+ 1\ep'}(w)X^+_{i\ep}(z)
\end{eqnarray*}
from which (\ref{E:R7}) follows in the case of $(\a_i|\a_j)=-1/2$. The case of
$(\a_i|\a_i)=1$ follows from (\ref{E:R11}).

For the case of $(\a_i|\a_j)=-1$, we have
\begin{eqnarray*}
&&(z-q^{-1}w)X^+_{n-1,\ep}(z)X^+_{n\ep'}(w)\\
&=&:X^+_{n-1,\ep}(z)X^+_{n\ep'}(w):
(z-q^{1+\ep'}_{\ep'}w)^{-1}
(z-q^{-1+\ep'}_{\ep'}w)^{-1}\\
&&\qquad\qquad\qquad\cdot (z-qw)(z-q^{-1}w)(z-w)\\
&=&(q^{-1}z-w)X^+_{n\ep'}(w)X^+_{n-1,\ep}(z),
\end{eqnarray*}
since it has no singularities.

The $-$ case of $(\a_n|\a_n)=2$ is shown as above, but
the $+$ case requires a different approach.
\begin{eqnarray*}
&&(z-q^2w)X^+_{n\ep}(z)X^+_{n\ep'}(w)\\
&=&:X^+_{n\ep}(z)X^+_{n\ep'}(w):q_{\ep}^{4\ep}
(z-q_{\ep}^{-1}q_{\ep'}q_{\ep}^{-\ep}q_{\ep'}^{\ep'}w)
(z-q_{\ep}q_{\ep'}^{-1}q_{\ep}^{-\ep}q_{\ep'}^{\ep'}w)\\
&&\quad\cdot (z-q_{\ep}^{-1}q_{\ep'}^{-1}q_{\ep}^{-\ep}q_{\ep'}^{\ep'}w)
(z-q_{\ep}q_{\ep'}q_{\ep}^{-\ep}q_{\ep'}^{\ep'}w)
(z-w)^{-1}(z-qw)^{-1}(z-q^{-1}w)^{-1}
\end{eqnarray*}
which implies that
\begin{eqnarray*}
&&\left((z-q^2w)X_n^+(z)X_n^+(w)-(w-q^2z)X_n^+(z)X_n^+(w)\right)
(q-q^{-1})^2(q^{1/2}-q^{-1/2})^2z^2w^2\\
&=&\sum_{\ep=\ep'}q^{\ep}:X^+_{n\ep}(z)X^+_{n\ep'}(w):
q^{2\ep}(z-w)\frac{(z-q_{\ep}^{2}w)(z-q_{\ep}^{-2}w)}
{(q-q^{-1})zw}\left(\delta(\frac{qw}z)-\delta(\frac{w}{qz})
\right)[2]_{q^{1/2}}^{2-2|\ep|}\\
&&\ \ +\sum_{\ep=-\ep'=\pm}:X^+_{n\ep}(z)X^+_{n\ep'}(w):
q^{2\ep}\frac{(z-q^{-\ep}w)(z-q^{-2\ep}w)}{z}\delta(q^{\ep}w/z)\\
&&\ \ -\sum_{\ep=0, \ep'=\pm}q^{\ep'/2}[2]_{q^{1/2}}
:X^+_{n\ep}(z)X^+_{n\ep'}(w):
\left((z-q^{2\ep'}w)+q^{2\ep'}(w-q^{-2\ep'}z)\right)\\
&&\ \ -\sum_{\ep=\pm, \ep'=0}q^{\ep/2}[2]_{q^{1/2}}
:X^+_{n\ep}(z)X^+_{n\ep'}(w):
\left((q^{2\ep}z-w)+(w-q^{2\ep}z)\right)\\
&=& \sum_{\ep\neq 0}:X^+_{n\ep}(z)X^+_{n,-\ep}(w):q^{2\ep}
\frac{(z-q^{-\ep}w)(z-q^{-2\ep}w)}{z}\delta(q^{\ep}w/z)\\
&&\ \ +:X^+_{n0}(z)X^+_{n0}(w):\frac{(z-w)(z-q^2w)(z-q^{-2}w)}{zw(q-q^{-1})zw}
\left(\delta(\frac{qw}z)-\delta(\frac{w}{qz})\right)[2]_{q^{1/2}}^2\\
&=&0
\end{eqnarray*}
since $:X^+_{n0}(z)X^+_{n0}(q^{\ep}z):=
:X^+_{n\ep}(z)X^+_{n,-\ep}(q^{\ep}w):$ for $\ep=\pm$.

The case of $X_n^-(z)X_n^-(w)$ is simpler than the "+" case, and is omitted.

{\it Proof of (\ref{E:R8}).} We only need to show the case when
$\a_i$ and $\a_j$ are connected in the Dynkim diagram.

\begin{eqnarray*}
&&(q^{-1}-q)(q^{1/2}-q^{-1/2})^2z^2w[X^+_{n}(z), X^-_{n-1}(w)]\\
&=&\sum_{\ep, \ep'}q^{\ep'/2}(-[2]_{q^{1/2}})^{1-|\ep|}\ep'
:X^+_{n\ep}(z)X^-_{n-1,\ep'}(w):\\
&&\ \ \left(\frac{(q_{\ep}^{1+\ep}z
-q^{\ep'/2}w)(q_{\ep}^{-1+\ep}z-q^{\ep'/2}w)}
{(z-q^{1/2}w)(z-q^{-1/2}w)}
-\frac{(q^{\ep'/2}w-q_{\ep}^{1+\ep}z)(q^{\ep'/2}w-q_{\ep}^{-1+\ep}z)}
{(w-q^{1/2}z)(w-q^{-1/2}z)}\right) \\
&=&\sum_{\ep\neq\ep'}q^{\ep'/2}(-[2]_{q^{1/2}})^{1-|\ep|}\ep'
:X^+_{n\ep}(z)X^-_{n-1,\ep'}(w):\\
&&\qquad \times\frac{(q^{\ep'/2}w-q_{\ep}^{1+\ep}z)
(q^{\ep'/2}w-q_{\ep}^{-1+\ep}z)}
{(q-q^{-1})zw}\left(\delta(\frac{q^{1/2}w}{z})-\delta(\frac{q^{-1/2}w}{z})
\right)\\
&=&
-\sum_{\ep=\pm}q^{\ep/2}\ep:X^+_{n\ep}(z)X^-_{n-1,-\ep}(q^{-\ep/2}z):
\frac{1-q^{-\ep}}{q^{-\ep}z}
\delta(\frac{q^{\ep/2}w}{z}) \\
&&\ -\sum_{\ep=\pm}
:X^+_{n0}(z)X^-_{n-1,\ep}(q^{-\ep/2}z):\frac{q^{-1/2}(1-q)}{q^{-\ep}z}
\delta(\frac{q^{\ep/2}w}{z})\\
&=&0
\end{eqnarray*}
since we have the following identities.
\begin{equation}
:X^+_{n\pm}(z)X^-_{n-1,\mp}(q^{\mp 1/2}z):
=:X^+_{n0}(z)X^-_{n-1,\pm}(q^{\mp 1/2}z):
\end{equation}

If follows from (\ref{E:R12}-\ref{E:R13}) that for $i=1, \cdots, n-1$
\begin{eqnarray*}
&&[X^+_{i}(z), X^-_{i}(w)]\\
&=&\frac{1}{(q^{1/2}-q^{-1/2})^2zw}
\sum_{\ep, \ep'}{\ep\ep'}:X^+_{i\ep}(z)X^-_{i\ep'}(w):
\left(\frac{z-q^{-\ep'/2}w}{q^{\ep/2}z-w}-\frac{w-q^{-\ep/2}z}
{q^{\ep'/2}w-z}\right)\\
&=&\frac{-1}{(q^{1/2}-q^{-1/2})^2zw}
\sum_{\ep=-\ep'}:X^+_{i\ep}(z)X^-_{i\ep'}(w):
\frac{z-q^{\ep/2}w}{w}\delta(\frac{q^{\ep/2}z}w)\\
&=&\frac{1}{(q^{1/2}-q^{-1/2})zw}\left(\psi_i(zq^{1/4})\delta(\frac{q^{1/2}z}w)
-\phi_i(zq^{-1/4})\delta(\frac{z}{q^{1/2}w})\right)
\end{eqnarray*}

Similarly using (\ref{E:R18}) and noting that $\ep=0$ does not contribute
to the commutator, we have
\begin{eqnarray*}
&&[X^+_{n}(z), X^-_{n}(w)]\\
&=&\frac{-1}{(q-q^{-1})^2z^2}
\sum_{\ep=\pm}q^{\ep/2}:X^+_{n\ep}(z)X^-_{n}(w):
\left(\frac{q^{-3\ep/2}(z-q^{3\ep/2}w)}{q^{\ep/2}z-w}-\frac{q^{-3\ep/2}
(q^{3\ep/2}w-z)}{w-q^{\ep/2}z}\right)\\
&=&\frac{-1}{(q-q^{-1})^2z^2}
\sum_{\ep=\pm}q^{\ep/2}:X^+_{n\ep}(z)X^-_{n}(w):
\frac{q^{-\ep}(z-q^{3\ep/2}w)}{w}\delta(\frac{q^{\ep/2}z}w)\\
&=&\frac{1}{(q-q^{-1})zw}\left(\psi_n(zq^{1/4})\delta(\frac{q^{1/2}z}w)
-\phi_n(zq^{-1/4})\delta(\frac{zq^{-1}}{w}\right)
\end{eqnarray*}

{\it Proof of Serre relations (\ref{E:R9}).}
Again we show the representative ones in the following.

\begin{eqnarray*}
X_{i\ep_1}^+(z_1)X_{i\ep_2}^+(z_2)X^+_{i+1,\ep}(w)&=&
:X_{i\ep_1}^+(z_1)X_{i\ep_2}^+(z_2)X^+_{i+1,\ep}(w):\\
&&\frac{(q^{\ep_1/2}z_1-q^{\ep_2/2}z_2)(z_1-q^{1/2}w)(z_2-q^{1/2}w)}
{(z_1-qz_2)(z_1-q^{\ep/2}w)(z_2-q^{\ep/2}w)}
\end{eqnarray*}

We have that
\begin{eqnarray}\label{E:S1}
&&X_{i\ep_1}^+(z_1)X_{i\ep_2}^+(z_2)X^+_{i+1,\ep}(w)-(q^{1/2}+q^{-1/2})
X_{i\ep_1}^+(z_1)X^+_{i+1,\ep}(w)X_{i\ep_2}^+(z_2)\nonumber\\
&&\ \ +X^+_{i+1,\ep}(w)X_{i\ep_1}^+(z_1)X_{i\ep_2}^+(z_2)\nonumber\\
&=&:X_{i\ep_1}^+(z_1)X_{i\ep_2}^+(z_2)X^+_{i+1,\ep}(w):
\frac{(q^{\ep_1/2}z_1-q^{\ep_2/2}z_2)}{(z_1-qz_2)
(z_1-q^{\ep/2}w)(z_2-q^{\ep/2}w)}\nonumber\\
&& \ \ \cdot \left((z_1-q^{1/2}w)(z_2-q^{1/2}w)
+[2]_{q^{1/2}}(z_1-q^{1/2}w)(w-q^{1/2}z_2)\right.\nonumber\\
&&\qquad\quad \left. +(w-q^{1/2}z_1)(w-q^{1/2}z_2)\right)\nonumber\\
&=&:X_{i\ep_1}^+(z_1)X_{i\ep_2}^+(z_2)X^+_{i+1,\ep}(w):
\frac{(q^{\ep_1/2}z_1-q^{\ep_2/2}z_2)(q^{-1/2}-q^{1/2})w}
{(z_1-q^{\ep/2}w)(z_2-q^{\ep/2}w)}
\end{eqnarray}
where we have used the identity
\begin{eqnarray}\label{E:S2}
(z_1-aw)(z_2-aw)
&+&(a+a^{-1})(z_1-aw)(w-az_2)+(w-az_1)(w-az_2) \nonumber\\
&=&(a^{-1}-a)w(z_1-a^2z_2)
\end{eqnarray}
Observe that the contracting factor in (\ref{E:S1}) is antisymmetric
under $(z_1, \ep_1)\mapsto (z_2, \ep_2)$, which implies the
Serre relation $X_i^+(z_1)X_i^+(z_2)X^+_{i+1}(w)+\cdots+(z_1
\leftrightarrow z_2)=0$.
Similarly we obtain that $X_i^+(z_1)X_i^+(z_2)X^+_{i-1}(w)-(q^{1/2}+q^{-1/2})
X_i^+(z_1)X^+_{i-1}(w)X_i^+(z_2)+\cdots=0$.

The Serre relation
\begin{eqnarray*}
X^-_{n}(z_1)X^-_{n}(z_2)X^-_{n-1}(w)&-&[2]_q
X^-_{n}(z_1)X^-_{n-1}(w)X^-_{n}(z_2)+X^-_{n-1}(w)X^-_{n}(z_1)X^-_{n}(z_2) \\
&& (z_1\leftrightarrow z_2)=0
\end{eqnarray*}
is proved similarly by considering
$X^-_{n}(z_1)X^-_{n}(z_2)X^-_{n-1, \ep}(w)$
and using the identity (\ref{E:S2}) with $a=q^{-1}$.

Finally let's check another representative Serre relation with $A_{ij}=-2$:
\begin{eqnarray*}
&&Sym_{z_1, z_2, z_3}
(X^+_{n-1}(z_1)X^+_{n-1}(z_2)X^+_{n-1}(z_3)X_{n}^+(w)-[3]_{q^{1/2}}
X^+_{n-1}(z_1)X^+_{n-1}(z_2)X_{n}^+(w)X^+_{n-1}(z_3)\\
&&\ +[3]_{q^{1/2}}
X^+_{n-1}(z_1)X_{n}^+(w)X^+_{n-1}(z_2)X^+_{n-1}(z_3)-
X_{n}^+(w)X^+_{n-1}(z_1)X^+_{n-1}(z_2)X^+_{n-1}(z_3))=0
\end{eqnarray*}

First we have
\begin{eqnarray*}
&&X^+_{n-1, \ep_1}(z_1)X^+_{n-1, \ep_2}(z_2)X^+_{n-1, \ep_3}(z_3)
X_{n, \ep}^+(w)\\
&=&:X^+_{n-1, \ep_1}(z_1)\cdots:\prod_{1\leq i<j\leq 3}
\frac{q^{\ep_i/2}z_i-q^{\ep_j/2}z_j}{z_i-qz_j}\prod_{i=1}^3
\frac{(z_i-qw)(z_i-w)}{(z_i-q_{\ep}^{1+\ep}w)(z_i-q_{\ep}^{-1+\ep}w)}
\end{eqnarray*}

Moving $X_{n, \ep}^+(w)$ around, we have that
\begin{eqnarray*}
&&X^+_{n-1, \ep_1}(z_1)X^+_{n-1, \ep_2}(z_2)X^+_{n-1, \ep_3}(z_3)
X_{n, \ep}^+(w)\\
&&\ \ -[3]_{q^{1/2}}X^+_{n-1, \ep_1}(z_1)X^+_{n-1, \ep_2}(z_2)
X_{n, \ep}^+(w)X^+_{n-1, \ep_3}(z_3)
\\
&&\ \ +[3]_{q^{1/2}}X^+_{n-1, \ep_1}(z_1)X_{n, \ep}^+(w)
X^+_{n-1, \ep_2}(z_2)X^+_{n-1, \ep_3}(z_3)\\
&&\ \ -X_{n, \ep}^+(w)X^+_{n-1, \ep_1}(z_1)X^+_{n-1, \ep_2}(z_2)
X^+_{n-1, \ep_3}(z_3)X_{n, \ep}^+(w)
\\
&=&:X^+_{n-1, \ep_1}(z_1)\cdots:\prod_{1\leq i<j\leq 3}
\frac{q^{\ep_i/2}z_i-q^{\ep_j/2}z_j}{z_i-qz_j}\prod_{i=1}^3
\frac{z_i-w}{(z_i-q_{\ep}^{1+\ep}w)(z_i-q_{\ep}^{-1+\ep}w)}\\
&&\ \ \cdot \left((z_1-qw)(z_2-qw)(z_3-qw)+[3]_{q^{1/2}}
(z_1-qw)(z_2-qw)(w-qz_3)\right.\\
&&\ \ \ \left.
[3]_{q^{1/2}}(z_1-qw)(w-qz_2)(w-qz_3)+(w-qz_1)(w-qz_2)(w-qz_3)\right)
\end{eqnarray*}
The expression in the parenthesis is simplified to be
$$
(q^{-1}-q)\left(w^2(z_1-(q+q^{-1})z_2+q^3z_3)+
w(z_1z_2-(q+q^{-1})z_1z_3+q^3z_2z_3)\right)
$$
Applying symmetrization of $S_3$ on $z_1, z_2, z_3$ (and on
$\ep_1, \ep_2, \ep_3$ accordingly) and factoring out the symmetric part,
we see that the proof of Serre
relation in the case is reduced to show the following identity
\begin{equation}\label{E:S3}
\sum_{\sigma\in S_3}sgn(\sigma) \sigma.(z_1-(q+q^2)z_2+q^3z_3)
\prod_{i<j}(qz_i-z_j)=0
\end{equation}
where the symmetric group $S_3$ acts on the ring of functions in $z_i$
($i=1, 2, 3$) in the natural way: $\sigma.z_i=z_{\sigma(i)}$.

The left-hand side (LHS) of (\ref{E:S3}) is a polynomial in $q$ of degree
5. It is easy to see that $q=\pm 1, 0, \infty$ are roots of the polynomial.
Let $\omega$ be a $3$rd primitive root of unity. We claim that
$\omega$ and $\omega^2$ are both roots of the LHS of (\ref{E:S3}).
In fact, plugging in $q=\omega$  reduces it into an
equivalent identity
\begin{eqnarray*}
&&\sum_{\sigma\in S_3}sgn(\sigma)\sigma.\prod_{i<j}(\omega z_i-z_j)\\
&=&\sum_{\sigma\in S_3}sgn(\sigma)\sigma.
\left((z_1^2z_2+z_1^2z+z_1z_2^2-z_2z_3^2)\right.\\
&&\ \ \left. +(1+2\omega)z_1z_2z_3+\omega(z_1^2z_3+z_1z_3^2+z_1z_2^2+z_2^2z_3)
\right)=0
\end{eqnarray*}
This proves the Serre relation and ends the proof of Theorem \ref{T:1}.

%==========================================================

\def\ch{{\rm{ch}}\,}
\def\Ker{{\rm{Ker}}\,}

   \section{Irreducible Representations}

   In this section, we will see the irreducible representations
   are realized in $\widetilde{\cal F}$.
   We treat the irreducible highest weight representations
   whose highest weights are the following four weights.
   $$ \mu_1 = - \frac{1}{2} \Lambda_0 , \;
      \mu_2 = - \frac{3}{2} \Lambda_0 + \Lambda_1 -\frac{\delta}{2}, \;
      \mu_3 = - \frac{1}{2} \Lambda_n+\frac{n\delta}8, \;
      \mu_4 = \Lambda_{n-1} -\frac{3}{2} \Lambda_n+\frac{n\delta}8. $$
   These weights are admissible
   and their characters are given by we the Weyl-Kac-Wakimoto character
formula \cite{kn:KW1}.
   In our case,
   the characters of the irreducible highest weight $\widehat{sp}_{2n}$-modules
   $L(\lambda)$ have the following forms \cite{kn:Lu}:
   $$ \ch L(\mu_1) + \ch L(\mu_2)
      = \frac{e^{-\frac{1}{2} \Lambda_0}}
             {\prod_{i=1}^{n}
              (p^{\frac{1}{2}} e^{ \vep_i};p)_{\infty}
              (p^{\frac{1}{2}} e^{-\vep_i};p)_{\infty}} \; , $$
   $$ \ch L(\mu_3) + \ch L(\mu_4)
      = \frac{e^{-\frac{1}{2} \Lambda_n} p^{-\frac{n}8}}
             {\prod_{i=1}^{n}
              (p e^{ \vep_i};p)_{\infty}
              (  e^{-\vep_i};p)_{\infty}} \; , $$
   where
   $p=e^{-\delta}$
   and
   $(a;p)_{\infty}=\prod^{\infty}_{n=0}(1-a p^n)$.

Define the operators $Q_i^-$ on $\widetilde{\cal F}^2_{\be}$
(the 2nd component)
by
\[ Q_i^-:=\oint Z^-_i(z) \frac{dz}{2 \pi \sqrt{-1}} .\]
We set subspaces ${\cal F}_i$ $(i=1,2,3,4)$
   of $\widetilde{\cal F}$ as follows.
\[\begin{array}{llll}
 {\cal F}_1 = \bigoplus_{\alpha \in Q}
                   {\cal F}'_{\alpha,\alpha},
    &&{\cal F}_2 = \bigoplus_{\alpha \in Q}
                   {\cal F}'_{\alpha+\vep_1,\alpha+\vep_1}\\
    {\cal F}_3 = \bigoplus_{\alpha \in Q}
                   {\cal F}'_{\alpha-\frac{1}{2}\lambda_n, \alpha},
   &&{\cal F}_4 = \bigoplus_{\alpha \in Q +\vep_n}
                   {\cal F}'_{\alpha-\frac{1}{2}\lambda_n, \alpha},
\end{array}\]
   where
   $$ {\cal F}'_{\alpha,\beta}
     ={\cal F}^{1}_{\alpha}
      \otimes
      \prod^n_{j=1} \Ker_{{\cal F}^{2}_{l_j \vep_j}} Q^-_j , $$
   for $ \beta=l_1 \vep_1 + \cdots + l_n \vep_n $.

   \bigskip
   \noindent
\begin{theorem}
   Each ${\cal F}_i$ $(i=1,2,3,4)$ is
   an irreducible highest weight $U_q$-module
   isomorphic to $V(\mu_i)$,
   The highest weight vectors are given by
   $ | \mu_1 \rangle
    =|0, 0 \rangle $,
   $ | \mu_2 \rangle
    = b_1(-1) | \lambda_1, \lambda_1 \rangle $,
   $ | \mu_3 \rangle
    =| -\frac{1}{2}\lambda_n, 0 \rangle $,
   $ | \mu_4 \rangle
    =|-\frac{1}{2}\lambda_n-\vep_n, -\vep_n \rangle $.
   \end{theorem}
   \bigskip
   \noindent
   {\it Proof.}
   Note that $[X^+_n(z), Q^-_n]=0$ is given in \cite{kn:K},
   we can see similarly that
   each $Q^-_i$ $(i=1, \cdots, n)$
   commutes or anticommutes with $X^{\pm}_j(z)$ for all $j=1, \cdots, n$.
   Therefore, each ${\cal F}_i$
   is a $U_q$-submodule of $\widetilde{\cal F}$.

   Next we calculate the character of ${\cal F}_i$.
   We may understand $Q^-_j$ as the zero mode $\eta_0$
   of the fermionic ghost system $(\xi, \eta)$ of dimension $(0,1)$.
   $$ \xi(z) =Z^+_j(z)=\sum_{n \in {\bf Z}} \xi_n z^{-n} \; , \quad
      \eta(z)=Z^-_j(z)=\sum_{n \in {\bf Z}}\eta_n z^{-n-1} . $$
   Since we have $\eta_0^2=0$ and $\xi_0 \eta_0 + \eta_0 \xi_0 = 1$,
   we obtain the following exact sequence:
   $$ 0 \longrightarrow \Ker_{{\cal F}^{2}_{l_j\vep_j}} Q_j^-
        \longrightarrow {\cal F}^{2}_{l_j\vep_j}
        \stackrel{Q_j^-}{\longrightarrow} {\cal F}^{2}_{(l_j-1)\vep_j}
        \stackrel{Q_j^-}{\longrightarrow} {\cal F}^{2}_{(l_j-2)\vep_j}
        \stackrel{Q_j^-}{\longrightarrow} \cdots . $$
   Using this exact sequence, we can compute the character of ${\cal F}_i$.
   Since $a_j(n)$ and $b_j(n)$ have weight $n\delta$, we get
   $$ \ch {\cal F}^{1}_{\alpha}
      =\frac{e^{\alpha} p^{-(\alpha|\alpha)}}
            {(p;p)_{\infty}^n} , \quad
     \ch {\cal F}^{2}_{l}
      =\frac{p^{\frac{1}{2}(l^2-l)}}
            {(p;p)_{\infty}} $$
   and hence
   \begin{eqnarray*}
   & & \ch({\cal F}_1 \oplus {\cal F}_2)=
   \ch(\bigoplus_{\alpha \in P} {\cal F}'_{\alpha,\alpha} ) \\
   &=& \ch(\bigoplus_{l_1, \cdots, l_n \in {\bf Z}}
           {\cal F}'
           _{l_1 \vep_1+\cdots+l_n \vep_n,
             l_1 \vep_1+\cdots+l_n \vep_n}) \\
   &=& \ch(\bigoplus_{l_1, \cdots, l_n \in {\bf Z}}
           {\cal F}^{1}_{l_1 \vep_1+\cdots+l_n \vep_n}
           \otimes \bigotimes^{n}_{j=1}
           \Ker_{{\cal F}^{2}_{l_j\vep_j}} Q^-_j) \\
   &=& \sum_{l_1, \cdots, l_n \in {\bf Z}}
       \left(
        \ch({\cal F}^{1}_{l_1 \vep_1+\cdots+l_n \vep_n})
        \prod^n_{j=1}
        \ch( \Ker_{{\cal F}^{2}_{l_j\vep_j}} Q^-_j)
       \right) \\
   &=& \sum_{l_1, \cdots, l_n \in {\bf Z}}
       \left(
        (p;p)_{\infty}^{-n} e^{l_1 \vep_1+\cdots+l_n \vep_n}
        \prod^{n}_{j=1}
        \left(
        (p;p)_{\infty}^{-1}
              \sum_{k \leq l_j}
               (-1)^{l_j-k} p^{-\frac{1}{2}(l_j^2-k^2+k)}
        \right)
       \right) \\
   &=& \sum_{l_1, \cdots, l_n \in {\bf Z}}
       \prod^{n}_{j=1}
       \left(
        (p;p)_{\infty}^{-2} e^{l_j \vep_j}
        \sum_{k \leq l_j}
        (-1)^{l_j-k} p^{-\frac{1}{2}(l_j^2-k^2+k)}
       \right) \\
   &=& \prod^{n}_{j=1}
       \left(
        (p;p)_{\infty}^{-2}
        \sum_{l \in {\bf Z}} e^{l \vep_j}
              \sum_{k \leq l}
               (-1)^{l-k} p^{-\frac{1}{2}(l^2-k^2+k)}
        \right) \\
   &=& \frac{e^{-\frac{1}{2} \Lambda_0}}
             {\prod_{i=1}^{n}
              (p^{\frac{1}{2}} e^{ \vep_i};p)_{\infty}
              (p^{\frac{1}{2}} e^{-\vep_i};p)_{\infty}} \; .
   \end{eqnarray*}
   The proof of the last equality was given in \cite{kn:K}.
   Similarly, we have
   $$ \ch({\cal F}_3 \oplus {\cal F}_4)
      = \frac{e^{-\frac{1}{2} \Lambda_l}}
             {\prod_{i=1}^{n}
              (p e^{ \vep_i};p)_{\infty}
              (  e^{-\vep_i};p)_{\infty}} \; , $$
   Hence, by the above character formulas,
   $$ \ch{\cal F}_i = \ch L(\mu_i) . $$
   As in Lusztig(\cite{kn:L}),
   ${\cal F}_i$ becomes a certain $\widehat{sp}_{2n}$-module
   in the classical limit $q \rightarrow 1$.
   Since the dimension of each weight space is invariant in the limit,
   the $\widehat{sp}_{2n}$-module in the classical limit is irreducible.
   Therefore ${\cal F}_i$ is irreducible for a generic $q$.

   Finally, it can be checked immediately
   that each  $| \mu_i \rangle$
   is a weight vector of weight $\mu_i$,
   and belongs to ${\cal F}_i$.
   $\mu_i$ is the highest weight of ${\cal F}_i$
   by the character.
   Hence $| \mu_i \rangle$ is a highest weight vector of ${\cal F}_i$.

   \hspace{\fill} $\Box$

   These representations become
   the ones constructed
   by Feingold and Frenkel \cite{kn:FF}
   in the classical limit $q \rightarrow 1$, since the latter is equivalent to
the $\beta-\gamma$ system of $\phi^1_j(z)$, $\phi^2_j(z)$ (cf \cite{kn:K} for
n=1 case).
   Our
   $a_j(n)$, $b_j(n)$ are related to them under $q\rightarrow 1$:
   $$ Y^{\pm}_j(z)
      \longrightarrow
      :e^{\pm 2 \phi^1_j(z)}: \; , \quad
      Z^{\pm}_j(z)
      \longrightarrow
      :e^{\pm   \phi^2_j(z)}: \; . $$

%======================================

\medskip

\centerline{Acknowledgements}
The first and third authors are supported in part by the NSA.
The second author is partly supported by
the Japan Society for Promotion of Science.

\end{document}